\documentclass{elsarticle}

\makeatletter
\def\ps@pprintTitle{%
 \let\@oddhead\@empty
 \let\@evenhead\@empty
 \def\@oddfoot{}%
 \let\@evenfoot\@oddfoot}
\makeatother
\usepackage{amsthm}
\usepackage{graphicx}
\usepackage{relsize}
\usepackage{exscale}
\usepackage{hyperref}
\usepackage{amsfonts}
\usepackage{epstopdf}
\usepackage{algorithmic}
\usepackage{amssymb,amsmath,epsfig,multirow}
\usepackage[displaymath, mathlines,running]{lineno}
\usepackage{mathtools} 
\usepackage{color}
\usepackage{lineno,hyperref}
\usepackage{rotating}
\usepackage{makeidx}
\usepackage{lscape}
\usepackage{float}
\usepackage{epsfig}
\usepackage{amsmath,amssymb}
\usepackage{enumitem}
\usepackage{epic,eepic}
\usepackage{overpic}
\usepackage{color}
\usepackage{amsfonts}
\usepackage{graphicx}
\usepackage{enumitem}
\usepackage{pgf,tikz,pgfplots}
\pgfplotsset{compat=1.15}
\usepackage{mathrsfs}

\usetikzlibrary{arrows}

\newtheorem{thm}{Theorem}
\numberwithin{thm}{section}

\newtheorem{cor}[thm]{Corollary}
\newtheorem{prop}[thm]{Proposition}

\newtheorem{defn}[thm]{Definition}
\newtheorem{remark}[thm]{Remark}

\newtheorem{assumption}[thm]{Assumption}

\newcommand{\del}{\delta}

\begin{document}
\begin{frontmatter}
\title{Approximation of Set-Valued Functions with images sets in $\mathbb{R}^d$}
\author{Nira Dyn}
\ead{niradyn@tauex.tau.ac.il}
\author{David Levin}
\ead{levindd@gmail.com}
\date{Received: date / Accepted: date}
\address{School of Mathematical Sciences. Tel-Aviv University. Tel-Aviv, Israel}

\begin{abstract}
{Given a finite number of samples of a continuous set-valued function F, mapping an interval to non-empty compact subsets of $\mathbb{R}^d$, $F: [a,b] \to K(\mathbb{R}^d)$, we discuss the problem of computing good approximations of F. We also discuss algorithms for a direct high-order evaluation of the graph of $F$, namely, the set
$Graph(F)=\{(t,y)\ | \ y\in F(t),\ t\in [a,b]\}\in K(\mathbb{R}^{d+1})$.
A set-valued function can be continuous and yet have points where the topology of the image sets changes. The main challenge in set-valued function approximation is to derive high-order approximations near these points.
In a previous paper, we presented with Q. Muzaffar, an algorithm for approximating set-valued functions with 1D sets ($d=1$) as images, achieving high approximation order near points of topology change. Here we build upon the results and algorithms in the $d=1$ case, first in more detail for the important case  $d=2$, and later for approximating set-valued functions and their graphs in higher dimensions.}
\end{abstract}


\end{frontmatter}

\section{Introduction}

In this paper we consider the approximation of set-valued functions from a finite collection of their samples. The class of set-valued functions we investigate consists of continuous (in the Hausdorff metric) set-valued functions, mapping $[a,b]\subset \mathbb{R}$ to $K(\mathbb{R}^d)$ - the collection of all nonempty compact subsets of $\mathbb{R}^d$.
Our goal is to present efficient algorithms for high-order approximation of $F$.
Recently, in \cite{MDL}, the authors and Q. Muzaffar, designed efficient algorithms for the interpolation of a set-valued function, mapping a closed interval to sets in $K(\mathbb R)$, given a finite number of samples of $F$,
which are 1D sets.
Such an algorithm computes a set-valued function that coincides with $F$ at the given samples. Analysis of the approximation error of the computed interpolants is carefully done as well.
The algorithms are inspired by the "metric polynomial interpolation" based on the theory in \cite{approximation_of_set_valued_functions}.     
 By this theory, a "metric polynomial interpolant" is a collection of polynomial interpolants to all the “metric chains” of the given samples of $F$. The algorithm
 in \cite{MDL} computes a small finite subset of "significant metric chains", which is sufficient for approximating $F$. More details of the algorithms in \cite{MDL}
are presented in the next section.
 
 In \cite{Levin1986}, the second author presented a method for set-valued interpolation for sets of general topology in any dimension $d$. The method is based on the interpolation of related signed-distance functions, which results in an implicit representation of the interpolant.
 A set-valued function can be continuous and yet have points where the topology of the image sets changes. A major challenge in set-valued function approximation is to derive high-order approximations near these points.
 The error analysis in \cite{Levin1986} is limited to closed sub-intervals in $[a,b]$ where the topology of $F$ does not change.
 In the simple case $d=1$ studied in \cite{MDL}, the authors present a comprehensive discussion on the global approximation error in $[a,b]$. One of the main contributions within \cite{MDL} lies in that it presents a method and an analysis for high-order approximation rate also near points of topology change. 

 In the present paper, we make use of the approach and the analysis in \cite{MDL} to achieve high-order approximation of set-valued functions with values
  in $K(\mathbb R^d)$. 
  We also use here the implicit function representation of a set, as suggested in \cite{Levin1986}.

  In the analysis of the approximation of real-valued functions, to achieve high-order approximation we usually assume that the function is sufficiently smooth. Here, for the approximation of a set-valued function $F$, we assume in this work that the graph of $F$, $Graph(F)=\{(t,y)\ | \ y\in F(t),\ t\in [a,b]\}\in K(\mathbb{R}^{d+1})$, has smooth boundaries. 

  An important notion for our analysis is the notion of points of change of topology in $Graph(F)$ (PCTs).
  We say the a point $p^*=(t^*,y^*)\in Graph(F)$ is a point of topology change of $F$ if there exists a neighborhood $U$ of $y^*$ in $\mathbb{R}^d$ and $\delta>0$ such only one of the two sets $F(t-\delta)\cap U$ and $F(t+\delta)\cap U$ is in $Graph(F)$.

In this work we consider, w.l.o.g, the approximation of a set-valued function $F: [0,1] \to K((0,1)^d)$ for $d\ge 1$
given its samples $\{F(ih)\}_{i=0}^N$, $h=1/N$.  
 We aim to generate an approximation $\tilde F(t)$ of the set $F(t)$ for any $t\in [0,1]$ and to achieve an error estimate of the form
$$d_{Haus}(\tilde F(t), F(t))=O(h^s), \text{\ as}\ h\to 0,$$
where $d_{Haus}$ denotes the Hausdorff distance in $K(\mathbb{R}^d) $, and $s > 1$ is the approximation order.

One way to represent the set $F(t)$ is in terms of its boundaries. The approximation of $F(t)$ in this representation requires the approximation of its boundaries. We use here the notion of the boundary of a compact set in $\mathbb{R}^n$ as
the set minus its interior.

A second way to represent F(t) is in terms of an {\bf inclusion algorithm} that determines whether a point in $\mathbb{R}^d$ belongs to the set $F(t)$ or not.

Here we use both approaches, i.e., we represent an approximation $\tilde F(t)$ to $F(t)$ by either computing its boundary or by an inclusion algorithm based upon an implicit representation. Namely, we construct a function $S:\mathbb{R}^d\to \mathbb{R}$, such that
\begin{equation}
\tilde F(t)=\{x\ :\ S(x)\ge 0\}.
\end{equation}

Another approximation target is the approximation of 
$Graph(F)$.
Here also the approximation may be by an explicit approximation of the boundaries of $Graph(F)$, or by an implicit inclusion algorithm. 

For $d\le 3$ there are known methods for solving our problem:
\begin{enumerate}

\item For $d=1$, the boundaries of $Graph(F)$ are curves in $\mathbb{R}^2$.
In \cite{MDL} an approximation order $O(h^s)$, as $h\to 0$, is obtained for $s>1$. The order $s$ depends upon the degree of the polynomials and the degree of splines used within the algorithms in \cite{MDL}, under the assumption that the boundaries of $Graph(F)$ are smooth enough.

\item For $d=2$, the boundaries of $Graph(F)$ are surfaces in $\mathbb{R}^3$. We remark that the case $d=2$ has many practical applications as it amounts to reproducing a 3D object (the graph of $F$) from its parallel 2D cross-sections
(the samples of $F$). A variety of methods devised for this case are in \cite{Levin1986}, \cite{LS}, \cite{GB}, \cite{KelsDyn}. None of these methods claim an approximation error higher than $O(h)$.

\item For $d=3$, the boundaries of $Graph(F)$ are 3D manifolds in ${\mathbb{R}^4}$ which are complicated to imagine. On the other hand, the evolution of $F$ as a function of $t$ is a familiar entity which is the animation of the changing 3D sets in between the given sets. At each frame of the animation, we display the boundary of the approximated 3D object $F(t)$. Based upon the method in \cite{Levin1986}, an efficient algorithm for solving this problem is presented in \cite{CSL}.

\end{enumerate}

\medskip

The approximation approach in this paper is built upon the significant progress achieved 
in \cite{MDL} for the approximation of set-valued functions with one dimensional sets as images.
Namely, the algorithm for approximating a function from $[0,1]$ to 
$K(\mathbb{R}^{d})$ is using multiple times an algorithm in \cite{MDL}.
In Section \ref{Pre} we review the algorithm presented in \cite{MDL} for the case $d=1$, which is the basic building block in our algorithm for the high dimensional approximation problem for $d\ge 2$. In particular, we focus on the high approximation order achieved.
This section contains also other preliminary algorithms used in 
constructing the approximations in high dimensions.
 In Section \ref{2Dcase}, we investigate in more detail the case $d=2$, which is related to the problem of reconstructing a 3D object from its 2D cross-sections. We show how the good approximation rates achieved in \cite{MDL} can be transferred to this case. In Section \ref{Highd} we extend our strategy for dealing with set-valued functions with images in higher dimensions.
We present an applicable procedure, based upon an implicit function representation, for computing high-order approximation of $d$-dimensional set-valued functions, $d\ge 3$. 

 In neighborhoods of points of topology change of the function, the approximation algorithm presented in \cite{Levin1986} fails to provide a high approximation order. On the other hand, the unique property of the 1D algorithms in \cite{MDL} is the ability to produce high-order approximations near points of topology change. By reducing the $d$-dimensional problem to a collection of one-dimensional problems, we achieve the desired approximation rate goals.

 \section{Preliminaries}\label{Pre}

 In this section, we bring known methods and results that are used in the construction and the analysis of the algorithms presented in the sequel.

\subsection{Review of the results for the case $d=1$.}\label{Pre1}
 Let us first present more details about the algorithms and the approximation results in \cite{MDL}, a paper that concerns the interpolation of a set-valued function $F : [a,b] \rightarrow K(\mathbb R^1)$ from a finite number of its samples.
 
It is assumed in \cite{MDL} that the graph of $F$ has a finite number of "holes". 


A (closed) hole of a set-valued function $F:[a,b]\to K(\mathbb{R})$ is a set of the form
    \begin{equation}
    \label{eq:definition_of_hole}
        H=\big\{(x,y):u(x)<y<v(x),x\in(c,d)\big\}\not\subset Graph(F),
    \end{equation}
    where $[c,d]\subset [a,b],\  u,v:[c,d]\to \mathbb{R}$, $u(c)=v(c)$, $u(d)=v(d)$ and   $u(x),v(x)\in F(x)$ for $x\in[c,d]$. We term the points $(c, u(c))$ and $(d, u(d))$ "points of change of topology of $F$ (PCT)" associated with the hole $H$.
The assumption that the number of holes in the graph of $F$ is finite, implies that each $F(t), t\in [a,b], $ is the union of a finite number of intervals. Note that a PCT is a point in the graph of $F$ such that the number of intervals of $F(t)$ to its left is different from this number to its right. 

The problem solved in \cite{MDL} is: Given a finite number of samples of a continuous set-valued function $F: [a,b] \to K(\mathbb R)$, namely $\{F(x_i)\}_{i=0}^N$, construct a set-valued function $\tilde{F}$ which interpolates the given samples of $F$, and derive the rate of approximation 
$$ d_{Haus}(F(t),\tilde{F}(t))=O(h^r),\ \text{\ as}\ h\to 0,$$
where $d_{Haus}$ denotes the Hausdorff metric in $K(\mathbb{R})$, $h=\max_{i=1}^N (x_i-x_{i-1})$ and $r>1$.

\medskip
The main steps of the approximation procedure in \cite{MDL} are as follows (w.l.o.g. we let $[a,b]=[0,1]$):
\begin{enumerate}
\item Detect points on the boundary of the graph of $F$ in the given samples of $F$.
\item Detect intervals $[x_i,x_{i+1}$ of topology change, i.e., such that the strip $[x_i,x_{i+1}]\times [0,1]$ contains a PCT of $F$.
\item Where there is no topology change, approximate the boundary by spline interpolation.
\item In a strip containing a PCT,
obtain high order approximation of the location of the PCT, for the following two cases:
\end{enumerate}

In Case A, it is assumed that the point of topology change is of type A, namely, the intersection of two smooth curves. Hence, this point is well approximated by the intersection of two polynomials which locally approximate the two curves.

In Case B, the point of topology change is of type B, namely, it is assumed that the boundary curves $u$ and $v$ are smooth and are H\"older 1/2 at the PCT. Namely, $u'$ and $v'$ are unbounded at the PCT. The location of the PCT is approximated as follows: First, we rotate the coordinate system by $\frac{\pi}{2}$, and then approximate the smooth boundary curve in the neighborhood of the PCT by fitting a local polynomial in the rotated system. The point of topology change is then approximated by the extrema of that polynomial.

In practice, a simple analysis of the data near a PCT can be used to distinguish between the two cases.

In the following, whenever referring to the methods in \cite{MDL} we refer to the algorithms and the analysis in Section 4 (Case A) and in Section 5 (Case B) therein. In both cases we use the topology change information to obtain a high-order approximation of the boundary close to the topology change location. More details can be found in \cite{MDL}.
\medskip

The following observations are important for our present work:
\begin{itemize}

\item{\it Approximation order}

In case A, assuming the boundary curves are $C^4$, an approximation order $O(h^4)$ is proved. This result can be easily extended to higher orders of approximation, $O(h^{k})$, under the assumption that the boundary curves are $C^{k}$.

In case B, assuming the boundary curves are $C^{2k}$, $k\ge 3$, an approximation order $O(h^{\frac{k}{2}})$ is achieved. 

\item{\it Reduced restrictions on $F$}

In \cite{MDL} it is assumed that the function $F$ is from $[0,1]$ into the class of {\it non-empty} subsets of $\mathbb{R}$. The left shape in Figure \ref{twoshapes} below describes the graph of such a function, with two holes, and four points of change of topology. The point marked as PCT1 is of type B, and PCT2 is of type A.
The function $F$ whose graph is shown in the right shape in Figure \ref{twoshapes},
depicts other types of points of change of topology, marked as PCT3 and PCT4. The function $F$ whose graph is on the right of Figure \ref{twoshapes} has another feature that is different from the functions discussed in \cite{MDL}. As $t$ varies from $t=0$ up to the location of PCT4, $F(t)$ is an empty set.
However, the restriction in \cite{MDL} to non-empty sets can be handled. The main issue is the detection of an interval containing a PCT. Then the same algorithm presented in \cite{MDL} for approximating the boundary of the graph near a PCT can be used for this case as well. Detecting the interval containing a PCT of the type PCT3 can be done by applying the algorithms in \cite{MDL} to the complements of the samples, namely to $\{F(t_i)^c\}$.

\begin{figure}[!ht]
  \centering
\begin{center}
{\includegraphics[width=1.1\textwidth]
{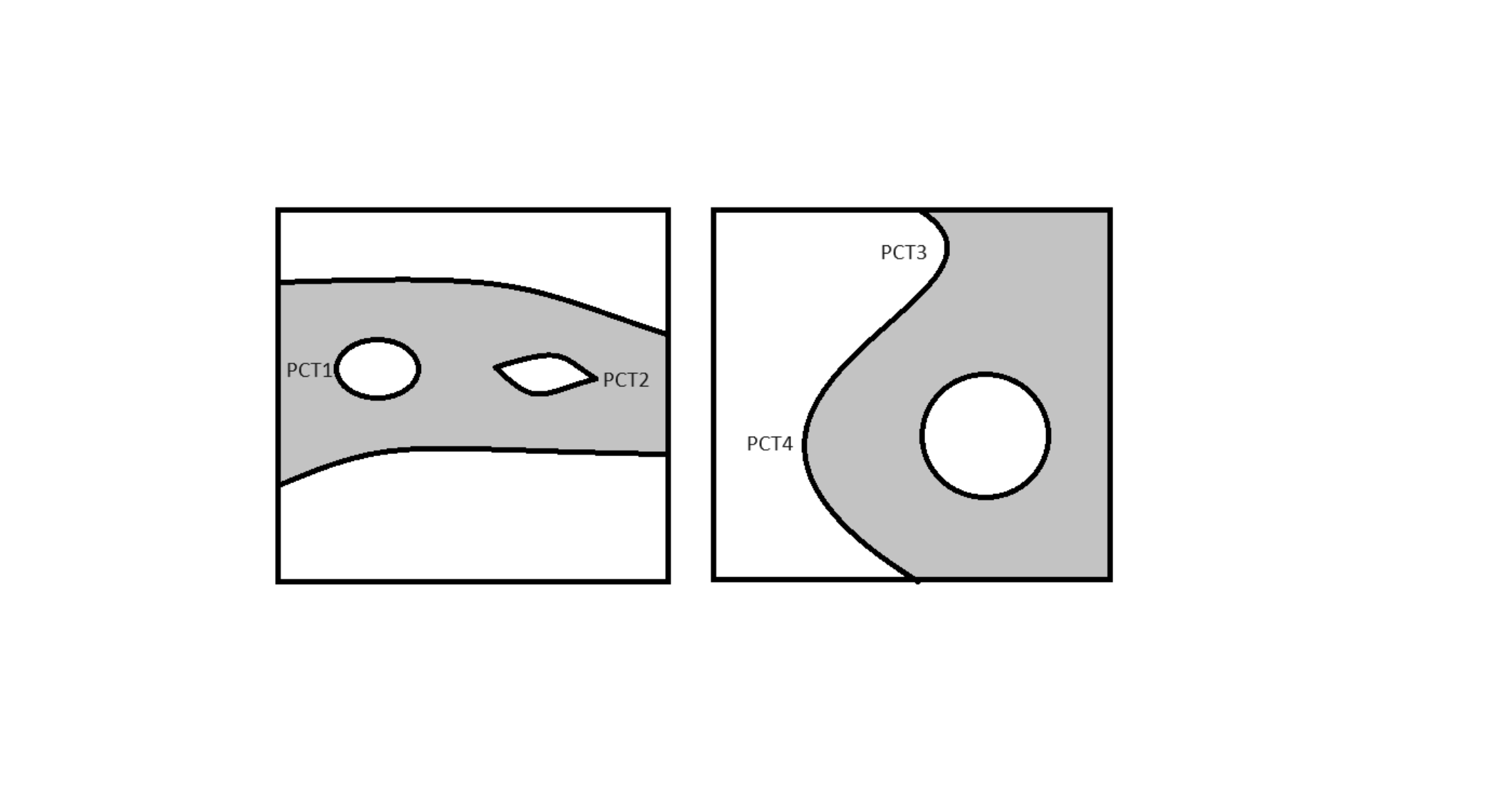}}\\
  \caption{The graphs of two examples of 1D set-valued functions, depicting different kinds of points of change of topology (PCT).}
  \label{twoshapes}
\end{center}
\end{figure}

\item{\it The approximation of boundary functions}

For later use, we point out that the approximation procedure in \cite{MDL} computes separate independent approximations to different parts of the boundary of $Graph(F)$.

\item{\it The approximation order $O(h^s)$}
\begin{defn}{\bf The approximation order $O(h^s)$.}\label{Ohs}

Each boundary curve computed by the algorithm in \cite{MDL} approximates the corresponding boundary curve within a certain approximation order.
Altogether, all the computed boundary curves provide an $O(h^s)$ approximation of the boundary of $Graph(F)$, where $s$ is the minimum of all the orders achieved in the approximations of the different parts of the boundary of $Graph(F)$.
\end{defn}

\end{itemize}

\subsection{Cross-sections}\label{CrossSection}

 Here we introduce the notion of cross-sections which is an essential tool for the approximation procedures we suggest in this paper.

\begin{figure}[!ht]
  \centering
\begin{center}
{\includegraphics[width=0.6\textwidth]
{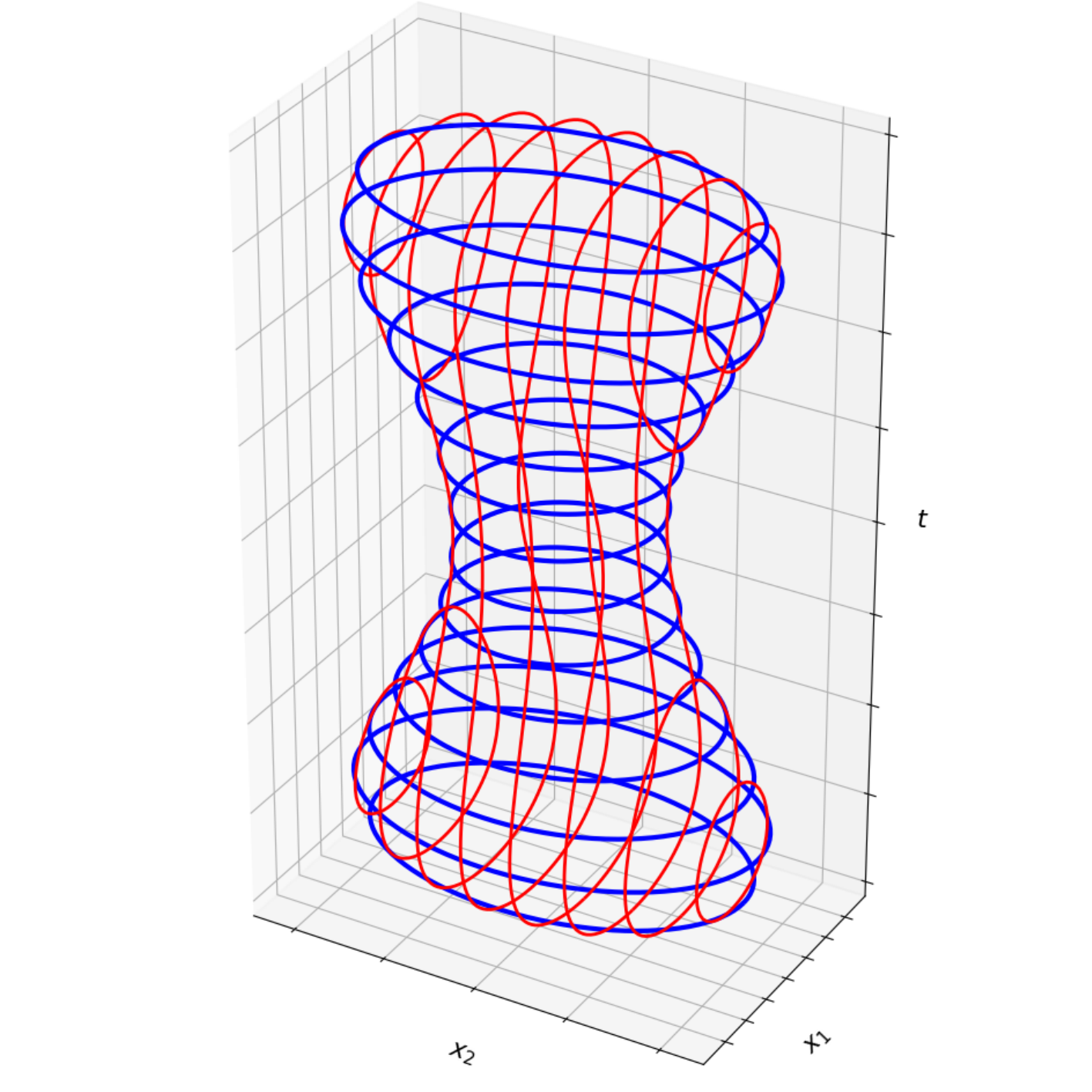}}\\
  \caption{A simple 3D object with cross-sections in two directions. The boundaries of the given cross-sections are in blue, and the boundaries of the cross-sections with constant $x_2$ are in red.}
  \label{cross-sections}
\end{center}
\end{figure}

\medskip
Let $F\ :\ [0,1]\to K(\mathbb{R}^d)$. 
We denote by 
$$F(t|x_k=\tau)\equiv F(t)\cap \{x_k=\tau\},$$
$\  \  \ \ $ the intersection of the $d$-dimensional set $F(t)$ with the hyperplane $x_k=\tau$.
In Figure \ref{cross-sections} we display the boundaries of cross-sections of $Graph(F)$ for $d=2$.

\bigskip

\medskip
Note that:
\begin{itemize}
\item For a fixed $\tau$, $G^{[k]}_\tau(t)\equiv F(t|x_k=\tau)$ is a set-valued function mapping $[0,1]$ to $K([0,1]^{d-1})$.

\item For a fixed $t$, $H^{[k]}_t(\tau)\equiv F(t|x_k=\tau)$ is another set-valued function mapping $[0,1]$ to $K([0,1]^{d-1})$.

Note that 
\begin{equation}\label{note1}
Graph(H^{[k]}_t)=F(t).
\end{equation}
\end{itemize}

\subsection{The sampling density condition}

\medskip
In the following, w.l.o.g., we assume that the samples of $F$ are taken at equidistant points, $\{t_i=ih\}_{i=0}^N$, $h=1/N$.

\begin{defn}\big[{\bf Condition on the sampling density $SD(h)$}\big]\label{SD(h)}

We assume that $h$ is 'small enough' so that we can capture the finest variations of $F$. The size 'small enough' $h$ is relative to the size $r_{min}$ of the minimal radius of curvature of the boundary surfaces of $Graph(F)$:
\begin{equation}\label{hcond}
h\le \sqrt{\epsilon r_{min}},
\end{equation}
where $\epsilon$ is small enough.
For short we say that it is required that condition $SD(h)$ is satisfied.
\end{defn}

Condition $SD(h)$ is related to the sampling density condition in the case of a piecewise linear approximation of a real-valued function $f$ on an interval, namely $h\le \sqrt{2\epsilon/max(f")}$, where $\epsilon$ is the required error bound.

\begin{defn}\big[{\bf Condition $SD1(h)$}\big]\label{DefnSD1(h)}

In \cite{MDL} an approximation order $O(h^s)$, as $h\to 0$, is proved, where $h$ is the sampling distance. The achievable order $s$ is of course limited by the smoothness order of the boundaries of $Graph(F)$. Then, the order $s$ depends upon the degree of the polynomials used to approximate the boundary near the neighborhood of points of topology change and upon the degree of splines used at regular regions. 

For example, let us view Step C of the approximation algorithm described in Section 5.2 of \cite{MDL}. 
It is assumed that the boundary curves are $C^{2k}$ smooth, and it is implicitly assumed that $h$ is small enough so that there are $2k$ sample points close to the point of topology change, for computing the polynomial $p_{2k-1}$ therein.
This condition is related to the more general condition $SD(h)$ and it can serve as a practical way of checking it for the case $d=1$. That is, near a point of topology change, we check whether there are enough sample points for applying the approximation algorithm in \cite{MDL}. We denote this condition as the $SD1(h)$ condition.
\end{defn}

\begin{remark}\label{Remarkmu}
Consider the 3D object shown in Figure \ref{cross-sections}. There are two planes $t=constant$ tangent to the object at its bottom and top points. As seen in the figure, within a significant neighborhood, $U_h$, of these points there is no data induced by the given samples of $F$. The size (diameter) of these neighborhoods depends upon the curvature of the surface of $Graph(F)$ at the tangency point.
If the principal curvatures are non-zero, the size of $U_h$ is $O(h^{0.5})$ as $h\to 0$, which corresponds to the diameter of a spherical cap section of height $h$ on the unit sphere.
$U_h$ is even larger if one of the principal curvatures is zero, namely, $O(h^\mu)$, $\mu<0.5$.

\end{remark}

\subsection{Spline quasi-interpolation operators}\label{QuasiInterpolation}

\medskip
In this paper, we use functions $S:\mathbb{R}^{d+1}\to \mathbb{R}$ for the implicit function representation of the 
 boundary surfaces of $Graph(F)$, $F:[0,1]\to K(\mathbb{R}^d)$, by spline-based quasi-interpolation operators applied to gridded data. The multivariate quasi-interpolation operators are obtained by tensor products of the univariate quasi-interpolation operators presented below.

Suppose $f\in\mathcal{C}^{p+1}(\mathbb{R})$,
$h>0$ a constant, $B_p(x)$ the $p$-degree B-spline supported on
$I_p =\left[-\frac{p+1}{2},\frac{p+1}{2}\right]$, with equidistant knots
$S_p = \left\{-\frac{p+1}{2},\hdots,\frac{p+1}{2}\right\},$ and define the vector
$f_{n,p}=(f_{n-\left\lfloor \frac{p}{2} \right\rfloor},\hdots,f_{n+\left\lfloor \frac{p}{2} \right\rfloor}),$
with $f_i=f(ih)$ and where $\left\lfloor \cdot \right\rfloor$ is the floor function. The following explicit local quasi-interpolation operator reproduces polynomials of degree $\le p$, and it approximates $f$ within order $O(h^{p+1})$ as $h\to 0$, see \cite{speleers}:
\begin{equation}\label{operatorQ}
Q_p(f)(x)=\sum_{n\in \mathbb{Z}}L_{p}(f_{n,p})B_p\left(\frac{x}{h}-n \right),
\end{equation}
where $L_p$ is the linear functional defined as:
\begin{equation}\label{operadorL}
L_p(f_{n,p})=\sum_{j=-\left\lfloor \frac{p}{2} \right\rfloor}^{\left\lfloor \frac{p}{2} \right\rfloor} c_{p,j} f_{n+j}.
\end{equation}
A general expression for the coefficients $c_{p,j}$, $j=-\left\lfloor \frac{p}{2} \right\rfloor,\hdots,\left\lfloor \frac{p}{2} \right\rfloor$ can be found in \cite{speleers}.

\subsection{Implicit approximation of smooth curves and surfaces}\label{Implicit}

\medskip

The subject of this work involves the approximation of sets with smooth boundaries. A constructive way of approximating a set is by approximating its boundaries, and an efficient approach, which is appropriate in any dimension, is to form an implicit approximation. That is, the approximation of a $d$-dimensional surface in $\mathbb{R}^{d+1}$, is defined as the zero-level surface of a function 
 $S$, $S:\mathbb{R}^{d+1} \to \mathbb{R}$. A valuable observation in this respect is the following curve approximation result appearing in \cite{ALR2}:
 
\begin{prop}\label{gammaappr1}
Assume $\Gamma$ is a $C^4$ smooth curve in $[0,1]^2$ with minimal curvature radius $>R$, and such that its $R$-neighborhood is not self-intersecting and not intersecting the boundaries of $[0,1]^2$. Let $P$ be a square mesh of points in  $[0,1]^2$,  of mesh size $\delta<R/10$. Let $\Gamma$ subdivide $[0,1]^2$ into two domains $\Omega_1$ and $\Omega_2$, and separate the mesh points $P$ into $P_1$ and $P_2$ respectively. Let us attach signed-distance values to the mesh points $P$ as follows:
$$D(q)=dist(q,\Gamma), \ \ \ q\in P_1,$$
$$D(q)=-dist(q,\Gamma), \ \ \ q\in P_2.$$
Consider $Q_3$ to be the quasi-interpolation operator from
vectors of values on $P$ to the space of bi-cubic splines on $[0,1]^2$.
Let $S_\epsilon$ be the bi-cubic spline defined by applying $Q_3$ to a perturbed signed-distance data
$$\tilde{D}(q)=D(q)+\epsilon_q,\ \ \ |\epsilon_q|\le \epsilon ,$$
and let $\Gamma_\epsilon$ be the zero-level curve of $S_\epsilon$. Then
\begin{equation}\label{disteps}
d_{Haus}(\Gamma,\Gamma_\epsilon)\le  C_0\delta^4+\|Q_3\|\epsilon.
\end{equation}
\end{prop}

The above result can be easily  extended to higher dimensions and higher approximation orders as follows:

\begin{prop}\label{gammaappr2}
Assume $\Gamma$ is a $C^m$ smooth surface in $[0,1]^d$ with minimal curvature radius $>R$, and such that its $R$-neighborhood is not self-intersecting and not intersecting the boundaries of $[0,1]^d$. Let $P$ be a square mesh of points in  $[0,1]^d$,  of mesh size $\delta<R/2m$. Let $\Gamma$ subdivide $[0,1]^d$ into two domains $\Omega_1$ and $\Omega_2$, and separate the mesh points $P$ into $P_1$ and $P_2$ respectively. Let us attach signed-distance values to the mesh points $P$ as follows:
$$D(q)=dist(q,\Gamma), \ \ \ q\in P_1,$$
$$D(q)=-dist(q,\Gamma), \ \ \ q\in P_2.$$
Consider $Q_{m-1}$ to be the quasi-interpolation operator from
vectors of values on $P$ to the space of $(m-1)$-order tensor product splines on $[0,1]^d$.
Let $S_\epsilon$ be the $(m-1)$-order tensor product spline defined by applying $Q_{m-1}$ to a perturbed signed-distance data
$$\tilde{D}(q)=D(q)+\epsilon_q,\ \ \ |\epsilon_q|\le \epsilon ,$$
and let $\Gamma_\epsilon$ be the zero-level surface of $S_\epsilon$. Then
\begin{equation}\label{disteps}
d_{Haus}(\Gamma,\Gamma_\epsilon)\le  C\delta^m+\|Q_{m-1}\|\epsilon.
\end{equation}
\end{prop}

\subsection{Estimating the distance of a point from a surface}\label{MLS}

\medskip
To compute the distances from a surface $\Gamma$ we suggest using the projection procedure presented in \cite{Levin2004} and extended to manifolds in \cite{MMLS}: 

Given points $X$ on an unknown surface $\Gamma$ and a point $p$ near the surface, the procedure finds a local reference plane and a projection of $p$ onto a point $T(p)$ near the surface. The projection is computed via a local least-squares by a polynomial patch $Q_p$. The distance $\|p-T(p)\|$ is suggested as a second-order approximation of $dist(p,\Gamma)$. A higher approximation order of $dist(p,\Gamma)$ may be defined as $dist(p,Q_p)$. The approximation order can be expressed in terms of the smoothness of $\Gamma$, the fill-distance of the set $X$, and the degree of the local polynomial approximation used. 

\begin{defn}{\bf Fill-distance}
    
The fill-distance of a set of points $X$ with respect to a surface (manifold) $\Gamma$ is the diameter of the biggest open ball centered at a point on $\Gamma$ that does not contain any point of $X$.

\end{defn}

\begin{prop}\label{MLS3}
Let $\Gamma$ be a $C^m$ $d$-dimensional surface, and let $X$ be a set of points on $\Gamma$ with fill-distance $kh$, where $k$ is a constant. Consider a point $p$ at distance $\le Kh$ from $\Gamma$, $K$ constant. Let $Q_p$ be the local polynomial patch defined by the projection procedure in \cite{MMLS}, of total degree $m-1$. Then,
\begin{equation}\label{distpGamma}
  dist(p,\Gamma)=dist(p,Q_p)+O(h^{m}),\ \text{as}\ h\to 0.
\end{equation}
\end{prop}

The proof follows by considering the approximation properties of local moving least-squares, as presented in \cite{MLS}.

\bigskip
\section{From 2D cross-sections to 3D objects}\label{2Dcase}

Let the set-valued function 
$F:[0,1]\rightarrow K([0,1]^2)$ be such that $Graph(F)$ has smooth boundaries. We let $t$ be the independent variable of $F$ and $x_1,x_2$ be the variables in $[0,1]^{2}$. Assuming we are given the 2D sets $\{F(ih)$, $i=0,...,N\}$, $h=1/N$, we would like to construct an approximation of $F(t)$ for any $t\in [0,1]$, and to approximate the set $Graph(F)$. Since $Graph(F)$ can be regarded as an object in $\mathbb{R}^3$, we consider here the problem of reconstructing a 3D object from its parallel 2D cross-sections.

 The assumption that $Graph(F)$ has smooth boundaries, implies that almost all 2D cross-sections of $Graph(F)$ also have smooth boundaries. We further assume that the radius of curvature of the boundary of $Graph(F)$ is bounded below by a constant $r^*_{min}$ and that $Graph(F)$ has a finite number of "holes". Hence also each 2D set $F(t)$ has a finite number of holes, and we can apply the results in \cite{MDL}. 
 
 The approximation procedure outlined below requires
 additional assumptions on $Graph(F)$ (Assumptions \ref{SDh1}, \ref{SDh2}).

\subsection{Approximating $F(t)$}\label{ApprF}

As declared in the introduction, the main challenge is to obtain high-order approximation globally, even near points of change of topology. We recall that the error analysis in \cite{Levin1986} is limited to closed sub-intervals in [0,1] where the topology of F does not change. An algorithm that overcomes the approximation problem near points of topology change has recently been presented in \cite{MDL} for functions $F:[0,1]\rightarrow K([0,1])$. To exploit the advanced results in \cite{MDL}, 
we consider the intersections of the sets $\{F(ih)\}_{i=0}^N$ with a line in the $x_1,x_2$ plane. W.l.o.g., we consider the one dimensional intersections of the sets $\{F(ih)\}_{i=0}^N$ with the lines $x_2=jh$, $j=0,...,N$. These are the 1D-sets 
$S_{i,j}=F(ih|x_2=jh)$, $i=0,...,N$, $j=0,...,N$. Each such  1D set is a union of a finite number of intervals.

Approximating $F(t)$ is a 2-stage procedure, where at each stage we use the 1D procedure in \cite{MDL}. In the first stage, we recover, with a high approximation rate, the cross-sections of the 3D object with the planes $x_2=jh$, $j=0,..., N$. Recalling the example in Figure \ref{cross-sections}, we recover the cross-sections depicted in red there. The first stage is required to obtain, for any $t$, new 1D data which can be used in the second stage to recover $F(t)$, using again the algorithm in \cite{MDL}.

The first stage in our procedure can be applied if the following assumption holds:
\begin{assumption}\label{SDh1}
We assume that the condition $SD1(h)$ (Definition \ref{DefnSD1(h)}) is satisfied for all the set-valued functions $F(\cdot|x_2=jh)$, $j=0,..,N$.
\end{assumption}
\noindent In case this assumption does not hold we suggest in Section \ref{Implicit} a different approach.

\medskip
{\bf First stage:} For a given $t$ and for a fixed $j$, we approximate the 1D set $F(t|x_2=jh)$, applying the 1D procedure developed in \cite{MDL} to the data $\{S_{i,j}\}_{i=0}^N$. The graph of the 1D set-valued function $F(\cdot|x_2=jh)$ is a 2D set. We recall that the procedure in \cite{MDL} finds an approximation to the boundaries of a 2D set from its 1D parallel cross-sections. Hence, following the approximation analysis in \cite{MDL}, the set-valued function $F(\cdot|x_2=jh)$ is recovered within $O(h^s)$ accuracy
(in the Hausdorff distance),
where $s$ can be chosen as explained in Definition \ref{DefnSD1(h)}. For the given $t$, and a fixed $j$, we obtain an approximation $\tilde F(t|x_2=jh)$ of the set of intervals $F(t|x_2=jh)$ with an approximation error $O(h^s)$.
We note that since the procedure in \cite{MDL} is interpolatory, then $$\tilde F(ih|x_2=jh)=F(ih|x_2=jh), \ i=0,...,N,\ j=0,...,N.$$

{\bf Second stage:} We consider the collection of 1D-sets $\{\tilde F(t|x_2=jh), j=0,..., N\}$ as another input for the 1D procedure in \cite{MDL}. Now we use it for approximating the boundary curves of the unknown 2D set $F(t)$. To analyze the error in the approximation of the boundary curves of $F(t)$, we should consider both the error within the 1D algorithm used in the second stage and the errors in the data $\tilde F(t|x_2=jh)$ due to the errors in the first stage. The approximated boundary curves define an approximation of the set $F(t)$. For a given value of $t$, to come up with an approximation error estimate as given in \cite{MDL}, we must assume also the following assumption:

\begin{assumption}\label{SDh2}
Condition $SD1(h)$ is satisfied for the set-valued function $H_t: [0,1]\to K([0,1])$,
$$H_t(\tau):=F(t|x_2=\tau).$$
\end{assumption}

\begin{prop} \label{Prophs}
Under Assumptions \ref{SDh1}
and \ref{SDh2}, the two-stage procedure above produces an approximation $\tilde F(t)\sim F(t)$ such that
\begin{equation}\label{dHausFt}
d_{Haus}(F(t),\tilde F(t))=O(h^s), \ \ {\text as\  } h\to 0.
\end{equation}
\end{prop}

\begin{proof}
By \cite{MDL}, using Assumption \ref{SDh1},
\begin{equation}\label{dHausFx2}
d_{Haus}(F(t|x_2=jh),\tilde F(t|x_2=jh))=O(h^s), \ \ {\text as\  } h\to 0.
\end{equation}
    Using the exact data $\{F(t|x_2=jh)\}_{j=0}^N$, under Assumption \ref{SDh2}, the procedure in \cite{MDL} produces an $O(h^s)$ approximation of the boundary curves of $F(t)$.
The data used in the second stage is an $O(h^s)$ perturbation of the exact data. Since the procedure in \cite{MDL} is linear w.r.t. the data, with bounded linear operators, an $O(h^s)$ perturbation of the data does not change the approximation order of the boundary curves of $F(t)$. This consequently implies the result (\ref{dHausFt}).

\end{proof}

\subsection{Approximating Graph(F) (the 3D object)}\label{3Dobject}

A practical approach for representing $Graph(F)$ is by displaying its boundaries. We discuss below four ways of approximating  $Graph(F)$ and representing the boundaries of $Graph(F)$:

\begin{enumerate}[label=\Roman*.]

\item {\bf Approximating the boundaries.}

Under Assumptions \ref{SDh1} and \ref{SDh2}, the above 2-stage process generates two groups of curves approximating the corresponding curves on the boundary of Graph(F). One group is related to cross sections of $Graph(F)$ orthogonal to the $x_2$-axis, namely the curves approximating the boundaries of $Graph(F(\cdot|x_2=jh))$, $j=0,..., N$. The other group of curves forms the boundaries of the given sets $F(ih)$, $i=0,..., N$.
By construction, the two groups of curves cross each other at the points that belong to $Boundary(S_{i,j})$, $i=0,..., N$, $j=0,..., M$. Thus, we have a net of curves lying at distance of order $O(h^s)$ from Boundary $(Graph(F))$.
A method for generating a smooth surface passing through a net of curves of general topology using a subdivision approach is presented in \cite{ALevin}. The analysis in \cite{ALevin} guarantees a $C^2$ smooth limit surface, however, it does not provide an estimate of the approximation error.

\item{\bf Approximating normals.}

A graphical representation of surfaces in 3D relies on the ability to compute reflections from the surface, i.e., the ability to compute the normal to the tangent plane at any point on the surface. Practically, it is enough to compute the normals at dense enough distributed points on the surface. The same way described above for approximating the boundary curves of cross-sections of the 3D body with planes of constant $x_2$, can be applied to construct approximations to the boundary curves of cross-sections of the 3D body with planes of constant $x_1$. This can be done for a dense enough collection of cross-sections and in an interpolatory manner. That is, the boundary curves of any cross-section in one direction interpolates the boundary curves of the given cross-sections $\{F(ih)\}_{i=1}^N$. At each intersection point of two curves, we can compute the tangent plane and the normal to the tangent plane using the tangents' directions of the two curves.

\item{\bf Using the metric average.}

We suggest the following procedure for obtaining an explicit construction of the
graph of a  Lipschitz continuous set-valued function $F: [0,1]\rightarrow K(\mathbb R^2)$ with error of order $O(h^s)$.
( $F$ is Lipschitz continuous if it satisfies
$$  d_{haus}(F(t_1),F(t_2)) \le L|t_1 - t_2|,\ \  t_1,t_2\in [0,1], $$
with $L$ a positive finite constant that depends on $F$). We further assume
that the boundaries of $Graph(F)$ are smooth except at points of change of topology.

In this procedure, we use the {\em metric average} between two sets in $K(\mathbb R^d)$ (see e.g. \cite{DFM}, Section 2.1). The metric average between the sets $A_1,A_2\in K(\mathbb R^d)$ with weight $w\in [0,1]$, is denoted by
$A_1 \oplus_w A_2 $, and is defined by 
$$A_1 \oplus_w A_2 = \{(1-w)a_1 + wa_2: (a_1,a_2)\in \Pi(A_1,A_2)\},$$
where $$\Pi(A_1,A_2)=\{(a_1,a_2)\in A_1\times A_2\ :\  \|{a_1 - a_2}\|\  \text {equals\ the \ distance\ from}\ a_1\  \text {to}\ A_2 \text{\ or\ the\   distance\ from}\  a_2\ \text{to}\  A_1\}. $$
The metric average has the "metric property" which is the property we use 
$$d_{Haus}(A_1 \oplus_w A_2,A_i)\le d_{Haus}(A_1,A_2),\ \ i=1,2.$$

 The first stage of this procedure is to compute approximations to the sets $F(\tau_j), \ \  j=0,\ldots,J$, with $\tau_{j+1}-\tau_j = h^s,\ \ j=0,\ldots\,J-1,\   \tau_0 = 0,\ \tau_J=1$. We denote the approximating set of $F(\tau_j)$ by ${\tilde F}(\tau_j)$ and compute it by the procedure for computing $F(t)$\ for $t\in[0,1]$ (see the previous section). For  the points where the data is given $t_0,\ldots,t_N$, no computation is needed.

\begin{prop}
Define for $t\in (\tau_j,\tau_{j+1})$,
$$ {\tilde F}(t)={\tilde F}(\tau_j)\oplus_w {\tilde F}(\tau_{j+1}),\ \ 
w=\frac{t-\tau_j}{h^s}\in (0,1).$$
Then, as $h\to 0$,
\begin{equation}\label{dHausF}
d_{Haus}(\tilde F(t),F(t))=O(h^s),\quad   t\in[0,1]
\end{equation}
and
\begin{equation}\label{dHausGraphF}
d_{haus}(Graph(\tilde F),Graph(F))=O(h^s).
\end{equation}
\end{prop}
\begin{proof}
First we prove (\ref{dHausF}).

For $t=\tau_j,\ j=0,\ldots,J$ the claim follows from the construction of the sets ${\tilde F}(\tau_j)$ in the first stage of our procedure, and from Proposition \ref{Prophs}. Thus,
\begin{equation}\label{dHausFtau}
d_{Haus}(\tilde F(\tau_j),F(\tau_j))=O(h^s),\ j=0,1,...,J.
\end{equation}

For $t\in (\tau_j,\tau_{j+1})$ for some $j\in \{0,1,2,\ldots, J-1\},$ the claim follows from the above, from the metric property of the metric average, and the Lipschitz continuity of $F$ using the inequality
$$
d_{Haus}({\tilde F}(t),F(t))\le d_{Haus}({\tilde F}(t),{\tilde F}(\tau_j)) + d_{Haus}({\tilde F}(\tau_j),F(\tau_j)) + d_{Haus}(F(\tau_j),F(t)).$$
 By the metric property of the metric average
 $d_{Haus}({\tilde F}(t),{\tilde F}(\tau_j))\le d_{Haus}({\tilde F} 
(\tau_{j+1}),{\tilde F}(\tau_j))$ and by the Lipschitz continuity of $F$,
$d_{Haus}(F(t),F(\tau_j)) \le L|t-\tau_j|=O(h^s)$. Using (\ref{dHausFtau}) we finaly get (\ref{dHausF}).

It remains to prove (\ref{dHausGraphF}).
By (\ref{dHausF}),
for any $t\in[0,1]$ there exist $p_t\in {\tilde F}(t)$ and $q_t\in F(t)$ such that
$$d_{Haus}({\tilde F}(t),F(t))=\|p_t-q_t\|= O(h^s). $$
Thus for any $t\in [0,1]$ the set $ \tilde {G}(t)=\{(t,y) : y\in {\tilde F}(t)\}\subset Graph({\tilde F})$ and the set
$G(t)=\{(t,y) : y\in F(t)\}\subset Graph(F)\}$ satisfy
$d_{Haus}(\tilde{G}(t),G(t))=\|(t,p_t)-(t,q_t)\|=\|p_t-q_t\|=O(h^s)$. Therefore \eqref{dHausGraphF} holds,
since,
$$d_{haus}(Graph(\tilde F),Graph(F))\le \max \{d_{Haus}(\tilde{G}(t),G(t)):\ t\in [0,1]\}.$$

\end{proof}

\item{\bf Implicit function representation.}

Another approximation strategy, that we use also for the high dimensional case, is based on the implicit function representation of a set \cite{Levin1986}, which yields an inclusion algorithm for the set.
First, we use the algorithm in \cite{MDL} to obtain approximations of boundary curves of cross-sections of the body (see Section \ref{CrossSection}) in two directions at a desired density. On these curves, we sample points at a desired density. These points are lying near the surface of the body at a distance of order $O(h^s)$ as $h\to 0$. Next, we use these points as input to the construction of an implicit representation of the surface. This is done by building a tensor product spline in 3D approximating the signed-distance function of the above-generated set of points. The spline function's zero-level surface would be the body surface's approximation, with an error $O(h^s)$. The detailed algorithm with the approximation error analysis is presented in Subsections \ref{Implicit}, \ref{Building} below.

\end{enumerate}

\begin{remark}\label{RemovingSDh}
 \big[Alleviating Condition $SD1(h)$\big]
 
It is possible to alleviate Assumptions \ref{SDh1} and \ref{SDh2}. We rely on the observation that the procedure in \cite{MDL} computes separate independent approximations to disconnected parts of the boundary curves of the graph of $F(\cdot|x_2=jh)$. Moreover, the $SD1(h)$ condition can be verified separately for each such part of the boundary curves, as explained in Definition \ref{DefnSD1(h)}.  Therefore, we compute those parts of approximating curves, for a given $x_2$, using the procedure in \cite{MDL}, only if Condition $SD1(h)$ is satisfied. The key observation here is that where a cross-section with a constant $x_2$ is close to being tangent to $Graph(F)$, some of the given cross-sections $\{F(ih)\}_{i=0}^N$  are close to being orthogonal to the boundary of $Graph(F)$. Thus, where a part of the boundary is missing due to the failure of the $SD1(h)$ condition for a curve of a cross-section with a constant $x_2$ (marked in red in the example in Figure \ref{cross-sections}), the missing information is obtained from the boundaries of the given samples $\{F(ih)\}$ (the blue curves in Figure \ref{cross-sections}).
\end{remark}

\subsection{The implicit function approach}\label{Implicit}

We recall that our main objective in this paper is to suggest an approximation method of high approximation order. Considering the alternative ways presented in 
Subsection \ref{3Dobject}, we choose to concentrate on the fourth way above involving the implicit function representation. There are three reasons to support this choice: First, an implicit function representation of a body is an exhaustive representation that contains all the information about the body and its boundaries. In particular, it enables approximating $F(t),\ t\in [0,1]$, without the assumption $SD1(h)$. Second, using this approach we can achieve the desired error estimates. Third, this approach works for any dimension, as presented in the next section. Therefore, we discuss this approach in more detail below.

For generating the implicit function representation we adopt the approach in \cite{Amat1}. In the first step, we need to distribute all around the boundaries of the body a set of points that are at a distance $O(h^s)$ from the boundary. These points should be densely distributed near the boundaries of the object to ensure a good approximation of the signed distance function. 

The main source of points are points on the boundary of the given samples. We define these points by distributing points along the boundaries of the given sets $\{F(ih)\}_{i=0}^N$ at distances of order $h$ from each other. We denote this collection of points by $Q^0_h$. We remark that this collection is not unique, and the way of finding it depends upon the representation of the given samples. For example, if the boundary curve parts of $F(ih)$ are given by explicit parametric representation, then points on the boundary can be obtained by evaluating the parametric representation at a sequence of parametric values. Another possibility is that the sets $\{F(ih)\}$ are given by an implicit function representation, namely,
$$ F(ih)=\{(x_1,x_2)\  |\  G_i(x_1,x_2)\ge 0,\ (x_1,x_2)\in [0,1]^2\},\ \  0\le i\le N.$$
In this case, we overlay a square net of points of mesh size $h$ over $[0,1]^2$. We search for vertices of the net that are in $F(ih)$ but their neighbors are not and use each of these vertices as a starting point for Newton's method for finding a zero of $G_i$.

By itself, $Q^0_h$ is not sufficient for defining a good enough approximation for the boundary of $Graph(F)$. These points do not cover well areas of the boundary $\Gamma$ of $Graph(F)$ near points where the boundary is tangent to the plane $t=constant$. This is shown well in the upper and the lower parts of the boundary of the object in Figure \ref{cross-sections}

To ensure complete coverage of the boundary of $Graph(F)$, we also compute points on the boundary curves of orthogonal cross-sections, namely the planes $x_2=jh$, $j=0,...,N$. By a 'complete coverage' we mean here that the collection of points has a fill-distance of order $O(h)$ with respect to $\Gamma$.

As described above in Section \ref{ApprF}, we can approximate the cross-section of the body with planes $x_2=jh$, and for each such cross-section, we obtain an approximation of its boundary curves. However, to ensure the desired error estimate, the condition $SD1(h)$ (see Definition \ref{DefnSD1(h)}) should be satisfied. Even for a small $h$, the condition $SD1(h)$ may not be satisfied for all the cross-sections of the body with the planes $x_2=jh$, $j=0,..., N$. This may occur near points where the plane $x_2=constant$ is nearly tangent to the surface of the body. Recalling Item 8 in Section \ref{Pre1}, the approximation algorithm in \cite{MDL} computes separate approximations for separate parts of the boundary curves in parametric form. For each of these parts, we check the $SD1(h)$ condition and discard those parts for which condition $SD1(h)$ fails.

On the computed boundary curves approximating the boundary of a cross-section $F(\cdot|x_2=jh)$, we sample points at distances $O(h)$. We denote this collection of points by $Q_{h,j}$. Recalling the points $Q_h^0$ distributed along the boundaries of the given samples, we define the augmented collection of points
$$Q_h^{[2]}=Q^0_h\mathsmaller{\bigcup}\cup_{j=0}^{N} Q_{h,j}.$$
These points lie at distance $O(h^s)$ from $\Gamma$, and their fill-distance with respect to $\Gamma$ is $O(h)$, as $h\to 0$.
We use this collection of points for defining the approximation of the signed-distance function. 

\subsection{Building and analyzing the signed-distance function approximation}\label{Building} 

We aim at constructing a tensor product spline $S(x,y)$, whose zero-level set defines the approximation to $\Gamma$. As in \cite{ALR2}, we look for $S$ which is an approximation to the signed-distance function from $\Gamma$. To construct $S$ we first overlay on $[0,1]^3$ a uniform 3D array of $(N+1)^3$ points $P$.

For each point $p\in P$ we approximate its distance from $\Gamma$ using the procedure described in Section \ref{MLS}.
Instead of data on the boundary of $\Gamma$, we have here the points $Q^{[2]}_h$ which lie at distance $O(h^s)$ from $\Gamma$.
Assuming a large enough $m$ in Proposition \ref{MLS3}, we obtain the approximated distances within error $O(h^s)$ as $h\to 0$.

Each point $p\in P$ lies on one of the given cross-sections $F(ih)$ of the 3D body, hence we know if it is inside or outside the body. To each point $p$ we assign the value of its approximate distance from $\Gamma$, with a plus sign if it belongs to the body and a negative sign if it is outside the body.

Using the quasi-interpolation operators presented in Section \ref{QuasiInterpolation}, a three-variate spline $S$ of degree $m$ is built with uniform knot's distance $h=1/N$, using the signed-distance values collected at the points of $P$. 

We denote the zero-level surface of the resulting $S$ as $\tilde\Gamma$, and this is the required approximation of $\Gamma$. The following approximation result follows directly from Proposition \ref{gammaappr2}.

\begin{cor}\label{gammacor}
Let the body's surface $\Gamma$ satisfy the smoothness assumptions as in Proposition \ref{gammaappr2} 
and let the approximated signed-distance values at the grid points $P$ be defined as above.
Let $S$ be the $m-th$-degree tensor product spline defined by quasi-interpolation applied to the approximated signed-distance values at $P$.  Denoting the zero-level surface of $S$ by $\tilde\Gamma$, we have
\begin{equation}\label{diste}
d_{Haus}(\Gamma,\tilde\Gamma)\le  C_1h^{m+1}+C_2h^s.
\end{equation}
\end{cor}

\section{Higher dimension procedures}\label{Highd}

We start with the example of approximating the graph of a set valued function, with values compact subsets of $[0,1]^d$ with $d=3$, since in this case the general pattern for any $d$ can be already observed.

\subsection{Computational algorithm for the case $d=3$}

In Section \ref{Implicit} we solve the approximation problem for the 2-dimensional problem, namely, approximating a function $F^{[2]}:[0,1]\to K([0,1]^2$ and its graph. The solution is based upon finding an implicit function representation for the boundaries of the 3D object which is the Graph of $F^{[2]}$.

Let us turn now to the 3-dimensional problem, namely, approximating the graph of a  function $F^{[3]}:[0,1]\to K([0,1]^3)$ given its samples $\{F^{[3]}(ih)\}_{i=0}^N$. Here as well, we would like to construct an implicit function approximation for the boundaries of $Graph(F^{[3]})$. To accomplish this we need to distribute all around the boundaries a collection of points of fill-distance $O(h)$. Given the approximation of the boundaries of $Graph(F^{[3]})$ as the zero-level set of $S(t,x_1,x_2,x_3)$, we can obtain the approximation $\tilde F^{[3]}(t)$ of $F^{[3]}(t)$ for any $t\in [0,1]$ as the set 
$$\tilde F^{[3]}(t)=\{(x_1,x_2,x_3) :  S(t,x_1,x_2,x_3)\ge 0\}.$$

Consider the cross-sections $F^{[2]}_{j_3}(\cdot)\equiv F^{[3]}(\cdot|x_3=j_3h)$, $0\le j_3\le N$. These are set-valued functions with images in $[0,1]^2$, and we are given their samples $F(ih|x_3=j_3h)$, $0\le i\le N$. Therefore, in principle, we can approximate these functions using the procedure in Section \ref{Implicit}. 
Recalling the approximation procedure there, it involves reconstructing cross-sections of $F^{[2]}$ using the 1D algorithm in \cite{MDL}. This is exactly what we suggest doing here, as explained below.

For each 2D cross-section $F_{j_3}^{[2]}$, we consider its 1D cross-sections
$$F^{[1]}_{j_2,j_3}(\cdot)\equiv F^{[2]}_{j_3}(\cdot|x_2=j_2h), \ 0\le j_2 \le N.$$
Given the samples of $F^{[3]}$, we can also find the samples of each function $F^{[1]}_{j_2,j_3}$. Now we can use the algorithm in \cite{MDL} to compute approximations to the boundary curves of $Graph(F^{[1]}_{j_2,j_3})$. As explained in Section \ref{Implicit}, we compute approximations only to those curves that satisfy the $SD1(h)$ condition. We apply this algorithm for each of the functions $F^{[1]}_{j_2,j_3}$, $0\le j_2,j_3 \le N$, 
obtaining approximations with accuracy of order $O(h^s)$ to their boundary curves. Along each of these approximations  to these boundary curves we distribute points at distance $O(h)$ from each other. All these points, together with the points $Q^0_h$ (defined in Section \ref{Implicit}), form the required set of points to define the signed-distance data, and the required tensor product spline function, $S:[0,1]^4\to \mathbb{R}$, for the implicit approximation of the boundaries of $Graph(F^{[3]})$. The full argument is presented in the next section for general dimension $d\ge 2$.

\medskip

We note that the 1D set-valued function $F^{[1]}_{j_2,j_3}$ can be written as
$$F^{[1]}_{j_2,j_3}(t)=F^{[3]}(t)\cap L_{j_2,j_3},$$
where $L_{j_2,j_3}$ is
the line segment in $\mathbb{R}^3$,
$$L_{j_2,j_3}\equiv \{(x_1,j_2h,j_3h): 0\le x_1 \le 1\}.$$
This observation is used in the next section to solve the approximation problem in higher dimensions.

\subsection{Computation and analysis  for higher dimensions}

Consider the $d$-dimensional problem, approximating a function $F^{[d]}:[0,1]\to K([0,1]^d)$ given its samples $\{F^{[d]}(ih)\}_{i=0}^N$.
We assume that the boundaries of $Graph(F^{[d]})$ are smooth, and $h$ is small enough such that the $SD(h)$ condition is satisfied.

As in the case $d=2$ (see subsection \ref {Implicit}), we distribute points on the boundaries of the given samples. We denote this collection of points by $Q^0_h$. Assuming the sets $\{F(ih)\}$ are given by an implicit function representation, namely,
$$ F(ih)=\{(x_1,x_2,...,x_d)\  |\  G_i(x_1,x_2,...,x_d)\ge 0,\ (x_1,x_2,...,x_d)\in [0,1]^d\},\ \  0\le i\le N,$$
we find these points as follows:
We overlay a mesh of points of mesh size $h$ over $[0,1]^d$. Then we search for vertices of the mesh that are in $F(ih)$ but some of their neighbors are not, and use each of these vertices as a starting point for Newton's method for finding a zero of $G_i$. The resulting points are distributed on the boundary of $F(ih)$ with fill-distance $O(h)$.

By itself, $Q^0_h$ is not sufficient for defining a good enough approximation for the boundary of $Graph(F)$. These points do not cover well areas of the boundary of $Graph(F)$ near points where the boundary is tangent to a plane $t=constant$. See, for example, the lower and the upper areas in Figure \ref{cross-sections}.

We define the collection of 1D set-valued functions
\begin{equation}\label{F1d}
F^{[1]}_{j_2,j_3,...,j_d}(t)=F^{[d]}(t)\cap L_{j_2,j_3,...,j_d},\ \ 0\le j_2,...,j_d \le N,
\end{equation}
where 
$$L_{j_2,j_3,...,j_d}\equiv \{(x_1,j_2h,j_3h,...,j_dh): 0\le x_1 \le 1\},$$
where $L_{j_2,j_3,...,j_d}$ is a line segment in $\mathbb{R}^d$.

Given the $d$-dimensional samples of $F^{[d]}$, we can find the 1D samples of each of the functions in (\ref{F1d}). Now we can use the algorithm in \cite{MDL} to compute approximations to the boundary curves of the graphs of the functions in (\ref{F1d}). As explained in Section \ref{Implicit}, we compute approximations only to those boundary curves that satisfy the $SD1(h)$ condition, thus obtaining
approximations to boundary curves with accuracy of order $O(h^s)$. Along each of these well-approximating curves, we distribute points at distance $O(h)$ from each other. We denote this collection of points by $Q_h^1\equiv \{Q_{h,j_2,j_3,...,j_d}$ $ 0\le j_2,...,j_d \le N\}$.

\begin{prop}
The collection of points

$$Q_h^{[d]}=Q^0_h\ \mathsmaller{\bigcup} \ Q_h^1,$$
is sufficient to define the signed-distance data, and the required spline function $S$ for the implicit approximation of the boundaries of $Graph(F^{[d]})$.
\end{prop}

\begin{proof}
Let $\Gamma$ denote the boundary of $Graph(F^{[d]})$.
The points $Q_h^0$, for small enough $h$, satisfy the $O(h)$ fill-distance property with respect to $\Gamma$ away from neighborhoods of points on $\Gamma$ where a hyperplane $t=constant$ is tangent to $\Gamma$.

Let $U_h\subset \Gamma$ be such a neighborhood. Due to the assumption that the boundary is smooth, the size of $U_h$ is of order $O(h^{\mu})$, $0<\mu\le 0.5$ (see Remark \ref{Remarkmu}).

For any choice of $j_2,j_3,...,j_d$, $0\le j_2,...,j_d \le N$, consider the 2D planar segment
$$H_{j_2,j_3,...,j_d}\equiv \{(t,x_1,j_2h,j_3h,...,j_dh): 0\le t,x_1 \le 1\}.$$
This planar segment is orthogonal to any hyperplane $t=constant$.
We observe that $Graph(F^{[1]}_{j_2,j_3,...,j_d})$ is contained in $H_{j_2,j_3,...,j_d}$. Therefore, for a choice of $j_2,j_3,...,j_d$, such that $H_{j_2,j_3,...,j_d}$ crosses $U_h$, $Graph(F^{[1]}_{j_2,j_3,...,j_d})$ is almost orthogonal to $U_h$. Hence, the radius of curvature of the boundary curves of $Graph(F^{[1]}_{j_2,j_3,...,j_d})$ crossing $U_h$ are bounded from below. To show this, we recall that $U_h$ is a neighborhood of a point $p$ on $Graph(F^{[d]})$ such that a hyperplane $t=constant$ is tangent to $\Gamma$ at $p$. Thus, $U_h$ may be described by a $d$-variate function over this  hyperplane. A boundary curve of $Graph(F^{[1]}_{j_2,j_3,...,j_d})$ crossing $U_h$ can be expressed as a univariate function, $u(x_1)$. Since $\Gamma$ is assumed to be smooth, $|u''(x)|$ is bounded, implying that the radius of curvature of such a boundary curve is bounded from below.

We conclude that for a small enough $h$, the boundary curves of $Graph(F^{[1]}_{j_2,j_3,...,j_d})$ crossing $U_h$ satisfy the $SD(h)$ condition, and therefore satisfies the $SD1(h)$ condition, and are thus computable and are of distance $O(h^s)$ from the boundary of $Graph(F^{[1]}_{j_2,j_3,...,j_d})$ and from $\Gamma$. 

Using all the functions in $(\ref{F1d})$, 
and distributing points at distances $O(h)$ along the approximating boundary curves of their graphs, we obtain a collection of points of fill-distance $O(h)$ distributed near all possible neighborhoods of the type of $U_h$.

Altogether, the points $Q^{[d]}_h$ form the required set of points for defining the signed-distance data, and the required spline function $S$ for the implicit approximation of the boundaries of $Graph(F^{[d]})$.
\end{proof}

Having the points $Q_h^{[d]}$ we construct a tensor product spline $S$, whose zero-level set defines the approximation of $\Gamma$. To construct $S$ we first overlay on $[0,1]^{d+1}$ a uniform grid of $(N+1)^{d+1}$ points $P$.
For each point $p\in P$ we approximate its distance from $\Gamma$ using the procedure described in Section \ref{MLS}.
Instead of data on the boundary of $\Gamma$ in Section \ref{MLS}, we use the points $Q^{[d]}_h$ which lie at distance $O(h^s)$ from $\Gamma$.
Assuming a large enough $m$ in Proposition \ref{MLS3}, we obtain the approximated distances within error $O(h^s)$ as $h\to 0$.

Each point $p\in P$ lies on one of the given cross-sections $F(ih)$, hence we know if it is inside or outside $Graph(F)$. To each point $p$ we assign the value of its approximate distance from $\Gamma$, with a plus sign if it is inside $Graph(F)$ and a negative sign otherwise.

Using the quasi-interpolation operators presented in Section \ref{QuasiInterpolation}, a $(d+1)$-variate spline $S$ of degree $m$ is built with uniform knots' distance $h=1/N$, using the signed-distance values collected at the points $P$. 

The approximation of $\Gamma$ is defined as the zero-level surface of $S$, and is denoted as $\tilde\Gamma$. The following approximation result follows directly from Proposition \ref{gammaappr2}.

\begin{thm}\label{gammacord}
Let the boundary surface $\Gamma$ of $Graph(F)$ satisfy the smoothness assumptions as in Proposition \ref{gammaappr2} 
and let the approximated signed-distance values at the grid points $P$ be defined as above.
Let $S$ be the tensor product spline of degree $m$ defined by quasi-interpolation applied to the approximated signed-distance values at $P$.  Denoting the zero-level surface of $S$ by $\tilde\Gamma$, we have
\begin{equation}\label{diste}
d_{Haus}(\Gamma,\tilde\Gamma)\le  C_1h^{m+1}+C_2h^s.
\end{equation}
\end{thm}

\end{document}